\documentclass[final,3p,11pt,nopreprintline]{elsarticle}

\usepackage{amsmath,amssymb}

\usepackage{amsthm}
\usepackage[T1]{fontenc}
\usepackage{lmodern}
\usepackage{etoolbox,needspace}
\usepackage[hidelinks]{hyperref}

\hypersetup{
  pdftitle={A co-sober Alexandrov space whose Scott power space is not co-sober},
  pdfauthor={Wei Ji and Xiaoquan Xu},
  pdfsubject={Co-sobriety and Scott power spaces},
  pdfkeywords={Co-sober space, compact saturated set, Scott topology, Smyth power poset, Alexandrov space}
}

\theoremstyle{plain}
\newtheorem{thm}{Theorem}[section]
\newtheorem{lem}[thm]{Lemma}
\newtheorem{prop}[thm]{Proposition}
\newtheorem{cor}[thm]{Corollary}
\theoremstyle{remark}
\newtheorem{rmk}[thm]{Remark}

\AtBeginEnvironment{thm}{\Needspace{5\baselineskip}}
\AtBeginEnvironment{lem}{\Needspace{5\baselineskip}}
\AtBeginEnvironment{prop}{\Needspace{5\baselineskip}}
\AtBeginEnvironment{cor}{\Needspace{5\baselineskip}}
\AtBeginEnvironment{rmk}{\Needspace{4\baselineskip}}

\newcommand{\Scott}{\sigma}
\newcommand{\Lower}{\omega}
\newcommand{\Alex}{\mathsf A}
\newcommand{\Q}{\mathcal Q}
\newcommand{\Fin}{\operatorname{Fin}}
\newcommand{\BoxOp}{\mathord\Box}
\newcommand{\N}{\mathbb N}

\begin{document}
\begin{frontmatter}

\title{A co-sober Alexandrov space whose Scott power space is not co-sober\tnoteref{funding}}
\tnotetext[funding]{This research was supported by the National Natural Science Foundation of China (Nos. 12471070, 12071199, 12561083).}

\author[guilin]{Wei Ji}
\ead{javeey@163.com}
\ead[url]{https://orcid.org/0000-0001-6220-7224}

\author[nanchang]{Xiaoquan Xu\corref{corresponding}}
\ead{xiqxu2002@163.com}
\ead[url]{https://orcid.org/0000-0003-1159-8477}
\cortext[corresponding]{Corresponding author: Xiaoquan Xu.}

\address[guilin]{School of Mathematics and Statistics, Guilin University of Technology, Guilin 541004, China}
\address[nanchang]{School of Computer Information Engineering, Nanchang Institute of Technology, Nanchang 330044, China}

\begin{abstract}
We show that the Scott power construction does not preserve co-sobriety. For a dcpo $P$, let $\operatorname{Fin}(P)$ be the poset of nonempty finitely generated upper sets ordered by reverse inclusion. We prove that co-sobriety of $\Sigma\operatorname{Fin}(P)$ implies co-sobriety of $\Sigma P$. The proof adapts the point-saturation and union argument used for full compact saturated power posets; the finite step follows from a result of Heckmann and Keimel. Applying this transfer to a dcpo with a non-co-sober Scott space yields a co-sober Alexandrov space whose Scott power space is not co-sober. We also construct a countable example for which the Smyth power poset is countable. In the countable example, the same poset is co-sober under the upper Vietoris topology but not under the Scott topology.
\end{abstract}


\end{frontmatter}

\section{Introduction}

In \cite[Definition~8.9 and Theorem~8.10]{ELS2004}, in order to provide an alternative approach to the theory of compactly generated spaces and establish a dual form of the Hofmann--Mislove Theorem, Escard\'o, Lawson and Simpson introduced the notion of co-sober spaces. A nonempty compact saturated subset is called $k$-irreducible if it is not the union of two proper compact saturated subsets. A $T_0$-space is co-sober if every $k$-irreducible compact saturated subset is the saturation of a point. Despite the formal duality with sobriety, the two notions behave differently; in particular, a sober space need not be co-sober~\cite[Example~2.9]{ShenWuXiZhao2020}; see also \cite[Example~2.8]{WenXu2018}.

\renewcommand{\thefootnote}{\fnsymbol{footnote}}

In an invited academic presentation\footnote[2]{X. Xu, Some open problems on co-sober spaces,
an invited academic presentation at the 11th International Symposium on Domain Theory and its Applications, Chengdu, China, 4--8 August 2026.}, Xu asked whether the Scott power space of a co-sober space is co-sober. In this paper, we give a negative answer to this question.

For a $T_0$ space $X$, let $\Q(X)$ be the set of its nonempty compact saturated subsets, ordered by reverse inclusion. The \emph{Scott power space} of $X$ is
\[
  P_{\Scott}(X)=\Sigma\Q(X),
\]
where $\Q(X)$ is endowed with its Scott topology. We study whether this construction preserves co-sobriety. The Scott topology on $\Q(X)$ can differ from the upper Vietoris topology, so the two power constructions must be treated separately.

Our first result shows that the Scott power construction does not preserve co-sobriety.

\begin{thm}\label{thm:main}
There exists a $T_0$ Alexandrov co-sober space $Y$ such that $P_{\Scott}(Y)$ is not co-sober.
\end{thm}

The construction can be made countable.

\begin{thm}\label{thm:countable-intro}
There exists a space $Y$ as in Theorem~\ref{thm:main} such that both $Y$ and $\Q(Y)$ are countable. For this example, the space $\Q(Y)$ endowed with the upper Vietoris topology is co-sober.
\end{thm}

The proof of Theorem~\ref{thm:main} uses a dcpo $P$ whose Scott space is sober but not co-sober, constructed by Shen, Wu, Xi and Zhao \cite[Example~2.9]{ShenWuXiZhao2020}. We show that non-co-sobriety of $\Sigma~\!\!P$ passes to $\Sigma\Fin(P)$, where
\[
  \Fin(P)=\{{\uparrow} F:\varnothing\ne F\subseteq P\text{ is finite}\}
\]
is ordered by reverse inclusion. If $Y=\Alex(P)$ is the Alexandrov space of the same poset, then $Y$ is co-sober and $\Q(Y)=\Fin(P)$.

Zhao, He and Wang proved related reflection results for full compact saturated power posets \cite[Theorem~3.13 and Corollary~3.16]{ZhaoHeWang2023}. Their argument uses point-saturation and union maps. We restrict these maps to the finite poset $\Fin(P)$. The additional input is
\[
  U\in\Scott(P)~\Longrightarrow~
  \BoxOp U=\{H\in\Fin(P):H\subseteq U\}\in\Scott(\Fin(P)),
\]
which follows from Heckmann and Keimel's finite Rudin result \cite[Corollary~3.9]{HeckmannKeimel2013}.

For the countable version, we start with a countable non-sober complete lattice, whose existence was proved by Miao, Xi, Li and Zhao \cite[Corollary~3.6]{MiaoXiLiZhao2023}. An argument using the lower topology supplies a nonprincipal co-sobriety witness. We then realize the lower space as a saturated subspace of the Scott space of a countable dcpo and apply the finite transfer.

Scott power spaces have also been studied in connection with sobriety and well-filteredness \cite{XuWenXi2023}. In a 2026 preprint, Xu, Xu, Ji and Qiu constructed a co-sober space whose upper Vietoris Smyth power space is not co-sober \cite{XuXuJiQiu2026}. Their result concerns a different topology.

\section{Preliminaries}

Let $P$ be a poset. For $A\subseteq P$, write
\[
  {\uparrow} A=\{y\in P:(\exists x\in A)\ x\le y\},
  ~
  {\downarrow} A=\{y\in P:(\exists x\in A)\ y\le x\}.
\]
We write ${\uparrow} x$ and ${\downarrow} x$ for the sets generated by singletons.

A subset $U\subseteq P$ is \emph{Scott open} if it is an upper set and, whenever a directed set $D\subseteq P$ has a supremum with $\bigvee D\in U$, one has $D\cap U\ne\varnothing$. This definition applies to arbitrary posets; no dcpo assumption is imposed. The Scott topology is denoted by $\Scott(P)$ and the resulting space by $\Sigma~\!\!P$. Its specialization order is the original order on $P$ \cite{GierzEtAl2003,GoubaultLarrecq2013}.

For a topological space $X$, a subset is saturated if it is an upper set in the specialization preorder. Throughout the paper, members of $\Q(X)$ are nonempty. A set $K\in\Q(X)$ is \emph{$k$-irreducible} if
\[
  K=K_1\cup K_2,
  ~ K_1,K_2\in\Q(X),
\]
implies $K=K_1$ or $K=K_2$. A $T_0$ space is \emph{co-sober} if each $k$-irreducible member of $\Q(X)$ is the saturation of a point \cite[Definition~8.9]{ELS2004}.

The Smyth order on $\Q(X)$ is
\[
  K\sqsubseteq L ~\Longleftrightarrow~ L\subseteq K.
\]
The upper Vietoris topology on $\Q(X)$ is generated by the sets
\[
  \BoxOp U=\{K\in\Q(X):K\subseteq U\},
  ~ U\in\mathcal O(X).
\]
We denote the corresponding space by $P_S(X)$. For a poset $P$, $\Alex(P)$ denotes the $T_0$ Alexandrov space whose open sets are all upper subsets of $P$.

The next two lemmas collect elementary facts needed below.

\begin{lem}\label{lem:alex-facts}
Let $P$ be a poset.
\begin{enumerate}
\item $\Q(\Alex(P))=\Fin(P)$ as Smyth ordered sets.
\item The space $\Alex(P)$ is co-sober.
\item The upper Vietoris topology on $\Fin(P)=\Q(\Alex(P))$ is the Alexandrov topology of the Smyth order. Hence $P_S(\Alex(P))$ is co-sober.
\end{enumerate}
\end{lem}

\begin{proof}
The saturated subsets of $\Alex(P)$ are the upper subsets of $P$. If an upper set $K$ is compact, then $\{{\uparrow} x:x\in K\}$ is an open cover of $K$, so $K={\uparrow} F$ for some nonempty finite $F\subseteq K$. Conversely, each principal upper set is compact: any open set containing its least element contains the whole principal upper set. A finite union of compact sets is compact. This proves (1).

If $K={\uparrow} x_1\cup\cdots\cup{\uparrow} x_n$ is $k$-irreducible, repeated application of the defining property gives $K={\uparrow} x_i$ for some $i$. This proves (2).

Every upper Vietoris open is upper in the Smyth order. Conversely, for $H\in\Fin(P)$ the set $H$ is open in $\Alex(P)$ and
\[
  \BoxOp H=\{K\in\Fin(P):K\subseteq H\}={\uparrow}_{\Fin(P)}H.
\]
Thus every principal Alexandrov open is upper Vietoris open, proving (3).
\end{proof}

\begin{lem}\label{lem:basic-transfer}
Let $X$ be a topological space.
\begin{enumerate}
\item The saturation of a compact subset of $X$ is compact.
\item If $A$ is a saturated subspace of a $T_0$ space $X$ and $K\subseteq A$ is compact saturated and $k$-irreducible in $A$, then $K$ has the same properties in $X$.
\end{enumerate}
\end{lem}

\begin{proof}
For (1), let an open cover of the saturation be given. A finite subfamily covers the original compact set. The union of this finite subfamily is open and contains the compact set, so it also contains its saturation.

For (2), compactness passes from the subspace to the ambient space. If $x\in K$ and $x\le_X y$, saturation of $A$ gives $y\in A$, and saturation of $K$ in $A$ gives $y\in K$. Thus $K$ is saturated in $X$. If $K=K_1\cup K_2$ with $K_i\in\Q(X)$, then $K_i\subseteq K\subseteq A$. Each $K_i$ is compact saturated in the subspace $A$, so $k$-irreducibility in $A$ gives $K=K_i$ for some $i$.
\end{proof}

\section{A finite transfer argument}

Let $P$ be a dcpo and order $\Fin(P)$ by reverse inclusion. The poset $\Fin(P)$ need not be a dcpo. Accordingly, every Scott topological argument below concerns only directed families whose supremum exists. The next proposition records the consequence of Heckmann and Keimel's Corollary~3.9 that is needed here.

\begin{prop}\label{prop:box-scott}
If $U$ is Scott open in $P$, then $\BoxOp U$ is Scott open in $\Fin(P)$.
\end{prop}

\begin{proof}
The set $\BoxOp U$ is upper in the Smyth order. Let $\mathcal H\subseteq\Fin(P)$ be directed, suppose that its supremum $H_*$ exists in $\Fin(P)$, and assume $H_*\subseteq U$. First observe that
\begin{equation}\label{eq:sup-intersection}
  H_*=\bigcap_{H\in\mathcal H}H.
\end{equation}
Indeed, $H_*\subseteq H$ for every $H\in\mathcal H$. Conversely, if $x$ belongs to the intersection, then ${\uparrow} x$ is an upper bound of $\mathcal H$ in the reverse inclusion order. The least upper bound property gives ${\uparrow} x\subseteq H_*$, hence $x\in H_*$.

For each $H\in\mathcal H$, choose a nonempty finite set $F_H$ with $H={\uparrow} F_H$. On nonempty finite subsets of $P$, write
\[
  F\preceq_{\mathrm{Sm}}G
  ~\Longleftrightarrow~
  {\uparrow} G\subseteq{\uparrow} F.
\]
Directedness of $\mathcal H$ under reverse inclusion implies that $\{F_H:H\in\mathcal H\}$ is directed for $\preceq_{\mathrm{Sm}}$. By \eqref{eq:sup-intersection},
\[
  \bigcap_{H\in\mathcal H}{\uparrow} F_H=H_*\subseteq U.
\]
Heckmann and Keimel's Corollary~3.9 \cite[Corollary~3.9]{HeckmannKeimel2013} gives an $H\in\mathcal H$ such that $F_H\subseteq U$. Since $U$ is upper, $H={\uparrow} F_H\subseteq U$. Therefore $\mathcal H\cap\BoxOp U\ne\varnothing$.
\end{proof}

Define
\[
  \eta:P\longrightarrow\Fin(P),~ \eta(x)={\uparrow} x.
\]
The map $\eta$ is monotone. If $D\subseteq P$ is directed, then $\eta[D]$ is directed and its supremum exists in $\Fin(P)$, with
\[
  \bigvee_{\Fin(P)}\eta[D]
  =\bigcap_{d\in D}{\uparrow} d
  ={\uparrow}\bigvee D.
\]
Thus $\eta$ preserves directed suprema and is Scott continuous as a map $\Sigma~\!\!P\to\Sigma\Fin(P)$.

For $K\in\Q(\Sigma~\!\!P)$ and $\mathcal A\subseteq\Fin(P)$, put
\[
  \Lambda(K)={\uparrow}_{\Fin(P)}\eta[K],
  ~
  \mu(\mathcal A)=\bigcup\mathcal A.
\]
These are the finite analogues of the point-saturation and union maps used for full compact saturated power posets; compare \cite[Theorem~3.13]{ZhaoHeWang2023}.

\begin{prop}\label{prop:lambda-mu}
Let $P$ be a dcpo.
\begin{enumerate}
\item If $K\in\Q(\Sigma~\!\!P)$, then $\Lambda(K)\in\Q(\Sigma\Fin(P))$ and $\mu(\Lambda(K))=K$.
\item If $\mathcal A\in\Q(\Sigma\Fin(P))$, then $\mu(\mathcal A)\in\Q(\Sigma~\!\!P)$.
\item For $K\in\Q(\Sigma~\!\!P)$, the set $K$ is $k$-irreducible if and only if $\Lambda(K)$ is $k$-irreducible.
\item For $K\in\Q(\Sigma~\!\!P)$, the set $K$ is principal if and only if $\Lambda(K)$ is principal in the Smyth order.
\end{enumerate}
\end{prop}

\begin{proof}
By Scott continuity of $\eta$, the image $\eta[K]$ is compact, and $\Lambda(K)$ is its saturation in $\Sigma\Fin(P)$. Lemma~\ref{lem:basic-transfer}(1) proves the first assertion in (1). If $H\in\Lambda(K)$, then $H\subseteq{\uparrow} x\subseteq K$ for some $x\in K$, so $\mu(\Lambda(K))\subseteq K$. Conversely, $x\in{\uparrow} x\in\Lambda(K)$ for every $x\in K$. Thus $\mu(\Lambda(K))=K$.

For (2), $\mu(\mathcal A)$ is nonempty and upper in $P$. Let $\{U_i:i\in I\}$ be a Scott open cover of $\mu(\mathcal A)$. Each $H={\uparrow} F\in\Fin(P)$ is a finite union of principal upper sets. Each principal upper set is compact in $\Sigma~\!\!P$, because every Scott open set containing its least element contains the whole principal upper set. Thus $H$ is compact. Therefore, for every $H\in\mathcal A$, there is a finite set $I_H\subseteq I$ such that
\[
  H\subseteq V_H:=\bigcup_{i\in I_H}U_i.
\]
By Proposition~\ref{prop:box-scott}, $\BoxOp V_H$ is Scott open in $\Fin(P)$. These sets cover $\mathcal A$. Compactness of $\mathcal A$ gives $H_1,\ldots,H_m\in\mathcal A$ with
\[
  \mathcal A\subseteq\BoxOp V_{H_1}\cup\cdots\cup\BoxOp V_{H_m}.
\]
It follows that $\mu(\mathcal A)$ is covered by the finitely many $U_i$ whose indices occur in $I_{H_1}\cup\cdots\cup I_{H_m}$. Thus $\mu(\mathcal A)$ is compact saturated.

Suppose that $K$ is $k$-irreducible and
\[
  \Lambda(K)=\mathcal A\cup\mathcal B,
  ~ \mathcal A,\mathcal B\in\Q(\Sigma\Fin(P)).
\]
Parts (1) and (2) give
\[
  K=\mu(\mathcal A)\cup\mu(\mathcal B).
\]
Assume $K=\mu(\mathcal A)$. For $x\in K$, choose $H_x\in\mathcal A$ with $x\in H_x$. Since $H_x$ is upper, ${\uparrow} x\subseteq H_x$, or $H_x\sqsubseteq{\uparrow} x$. Saturation of $\mathcal A$ gives ${\uparrow} x\in\mathcal A$. Hence $\eta[K]\subseteq\mathcal A$, and another use of saturation gives $\Lambda(K)\subseteq\mathcal A$. Therefore $\Lambda(K)=\mathcal A$.

Conversely, suppose that $\Lambda(K)$ is $k$-irreducible and $K=K_1\cup K_2$ with $K_i\in\Q(\Sigma~\!\!P)$. Then
\[
  \Lambda(K)=\Lambda(K_1)\cup\Lambda(K_2).
\]
Thus $\Lambda(K)=\Lambda(K_i)$ for some $i$, and applying $\mu$ gives $K=K_i$. This proves (3).

If $K={\uparrow} x$, then
\[
  \Lambda(K)={\uparrow}_{\Fin(P)}\eta(x).
\]
Conversely, suppose $\Lambda(K)={\uparrow}_{\Fin(P)}H_0$. Every member of ${\uparrow}_{\Fin(P)}H_0$ is contained in $H_0$, and $H_0$ itself belongs to this set. Hence its union is $H_0$, and (1) gives $K=H_0$. Since $H_0\in\Lambda(K)$, there is $x\in K$ with $H_0\subseteq{\uparrow} x$. Because $x\in H_0$ and $H_0$ is upper, ${\uparrow} x\subseteq H_0$. Hence $K={\uparrow} x$. This proves (4).
\end{proof}

\begin{cor}\label{cor:finite-reflection}
For every dcpo $P$,
\[
  \Sigma\Fin(P)\text{ co-sober}
  ~\Longrightarrow~
  \Sigma~\!\!P\text{ co-sober}.
\]
Equivalently, non-co-sobriety of $\Sigma~\!\!P$ implies non-co-sobriety of $\Sigma\Fin(P)$.
\end{cor}

\begin{proof}
If $\Sigma~\!\!P$ is not co-sober, choose a nonprincipal $k$-irreducible member $K$ of $\Q(\Sigma~\!\!P)$. Proposition~\ref{prop:lambda-mu} shows that $\Lambda(K)$ is a nonprincipal $k$-irreducible member of $\Q(\Sigma\Fin(P))$.
\end{proof}

\begin{rmk}
The proof does not assert that $\Lambda$ is Scott continuous as a map between compact saturated power posets, or that $\Lambda$ and $\mu$ form a topological retraction. Proposition~\ref{prop:box-scott} is used to prove compactness of $\mu(\mathcal A)$ and supplies the finite step not contained in the identity $\mu\Lambda=\mathrm{id}$.
\end{rmk}

\section{A counterexample}

\begin{proof}[Proof of Theorem~\ref{thm:main}]
Shen, Wu, Xi and Zhao constructed a dcpo $P$ such that $\Sigma~\!\!P$ is sober but not co-sober \cite[Example~2.9]{ShenWuXiZhao2020}. Let
\[
  Y=\Alex(P).
\]
By Lemma~\ref{lem:alex-facts}, $Y$ is a $T_0$ co-sober space and $\Q(Y)=\Fin(P)$. Corollary~\ref{cor:finite-reflection} shows that
\[
  P_{\Scott}(Y)=\Sigma\Q(Y)=\Sigma\Fin(P)
\]
is not co-sober.
\end{proof}

The preceding argument uses no cardinality property of the dcpo in \cite[Example~2.9]{ShenWuXiZhao2020}. The remainder of the paper gives a separate countable construction.

\section{A countable witness in a lower topology}

Miao, Xi, Li and Zhao proved that a countable non-sober complete lattice exists \cite[Corollary~3.6]{MiaoXiLiZhao2023}. Fix such a lattice $L$. Since $\Sigma~\!\!L$ is $T_0$ and non-sober, there is a nonempty irreducible Scott closed set $C\subseteq L$ such that
\begin{equation}\label{eq:C-nonprincipal}
  C\ne{\downarrow} x~\text{for every }x\in L.
\end{equation}
Only the existence statement from \cite[Corollary~3.6]{MiaoXiLiZhao2023} is used; no witness from that paper is imported.

The lower topology $\Lower(L)$ is generated by the subbasic sets
\[
  L\setminus{\uparrow} a,~ a\in L.
\]
Its specialization order is the reverse of the lattice order. Thus its saturated sets are precisely the lower subsets of $L$, and the saturation of $x$ is ${\downarrow} x$.

We shall also use the elementary inclusion
\begin{equation}\label{eq:lower-scott-op}
  \Lower(L)\subseteq\Scott(L^{\mathrm{op}}).
\end{equation}
Indeed, each lower subbasic set $L\setminus{\uparrow} a$ is an upper set in $L^{\mathrm{op}}$. If a filtered set $F\subseteq L$ has meet $\bigwedge F\notin{\uparrow} a$, then some $f\in F$ satisfies $a\nleq f$; otherwise $a$ would be a lower bound of $F$ and hence $a\leq\bigwedge F$. Thus $L\setminus{\uparrow} a$ is Scott open in $L^{\mathrm{op}}$, and so is every lower open set.

\begin{lem}\label{lem:lower-compact}
Let $L$ be a complete lattice and let $K\subseteq L$ be nonempty. Then $K$ is compact saturated in $(L,\Lower(L))$ if and only if $K$ is Scott closed in $L$.
\end{lem}

\begin{proof}
Suppose first that $K$ is Scott closed. It is lower and therefore saturated. By Alexander's subbase theorem, it is enough to consider a cover by lower subbasic opens. Let
\[
  K\subseteq\bigcup_{i\in I}(L\setminus{\uparrow} a_i),
\]
and suppose that no finite subfamily covers $K$. For every finite $F\subseteq I$, choose $k_F\in K$ with $a_i\le k_F$ for all $i\in F$, and put
\[
  d_F=\bigvee_{i\in F}a_i,
\]
where $d_\varnothing=\bot$. Since $d_F\le k_F$ and $K$ is lower, $d_F\in K$. The family $\{d_F:F\subseteq I\text{ finite}\}$ is directed, and Scott closedness gives
\[
  \bigvee_{\substack{F\subseteq I\\F\text{ finite}}}d_F
  =\bigvee_{i\in I}a_i\in K.
\]
This element lies in every ${\uparrow} a_i$, contradicting the cover. Hence $K$ is compact.

Conversely, suppose that $K$ is compact saturated. Then $K$ is lower. Let $D\subseteq K$ be directed and set $s=\bigvee D$. If $s\notin K$, then for every $k\in K$ there is $d\in D$ with $d\nleq k$; otherwise $k$ would be an upper bound of $D$, and lowerness would give $s\in K$. Hence
\[
  \{L\setminus{\uparrow} d:d\in D\}
\]
is an open cover of $K$. A finite subcover corresponds to $d_1,\ldots,d_n\in D$. Choose $d\in D$ above all $d_i$. Then $d\in K$ belongs to none of the selected complements, a contradiction. Thus $s\in K$, and $K$ is Scott closed.
\end{proof}

A nonempty closed set is irreducible exactly when it is not the union of two proper closed subsets. Lemma~\ref{lem:lower-compact} and \eqref{eq:C-nonprincipal} give the following.

\begin{cor}\label{cor:lower-witness}
The set $C$ is a nonprincipal $k$-irreducible compact saturated subset of $(L,\Lower(L))$.
\end{cor}

\section{A countable staircase dcpo}

Fix a sequence $a_0,a_1,\ldots$ in which every element of $L$ occurs. For $x\in L$ and $n\in\N$, define
\[
  U_n(x)=
  \bigcap_{\substack{i\le n\\a_i\nleq x}}
  (L\setminus{\uparrow} a_i).
\]
The empty intersection is $L$. Each $U_n(x)$ is a lower open neighbourhood of $x$, the sequence is decreasing, and
\[
  \bigcap_{n\in\N}U_n(x)={\downarrow} x.
\]
Moreover, the sets $U_n(x)$ form a neighbourhood base at $x$ in the lower topology, since every basic lower open neighbourhood involves only finitely many subbasic complements.

Let
\[
  D_L=\{\widehat x:x\in L\}\mathbin{\dot\cup}
      \{x_n:x\in L,\ n\in\N\}.
\]
Define a relation on $D_L$ by
\begin{align*}
  x_n\le x_m
  &~\Longleftrightarrow~ n\le m
  ~\text{for the same label }x, \\
  x_n\le\widehat y
  &~\Longleftrightarrow~ y\in U_n(x), \\
  \widehat x\le\widehat y
  &~\Longleftrightarrow~ y\le_L x.
\end{align*}
We impose no other comparisons. Reflexivity and antisymmetry are immediate. Transitivity in the mixed cases follows from $U_m(x)\subseteq U_n(x)$ for $n\le m$ and from lowerness of $U_n(x)$. Thus this relation is a partial order, and $D_L$ is countable.

\begin{lem}\label{lem:DL-dcpo}
The poset $D_L$ is a dcpo. Let $E\subseteq D_L$ be directed.
\begin{enumerate}
\item If $E$ contains a point $\widehat x$, set
\[
  E^\sharp=\{x\in L:\widehat x\in E\}.
\]
Then $E^\sharp$ is filtered and
\[
  \bigvee E=\widehat{\bigwedge E^\sharp}.
\]
\item If $E$ contains no point of the form $\widehat x$, then all its elements lie on one staircase. If the occurring indices have a maximum $m$, the supremum is $x_m$; if they are unbounded, the supremum is $\widehat x$.
\end{enumerate}
\end{lem}

\begin{proof}
Suppose that $E$ contains $\widehat{x_0}$. For each $e\in E$, directedness gives $u\in E$ above both $e$ and $\widehat{x_0}$. No staircase point lies above a point $\widehat x$, so the point part of $E$ is cofinal. If $x,y\in E^\sharp$, a common upper bound in $E$ has the form $\widehat z$ with $z\le_L x,y$. Hence $E^\sharp$ is filtered. If $m=\bigwedge E^\sharp$, then $\widehat m$ is the least upper bound of the cofinal point part and therefore of $E$.

Now suppose that $E$ contains only staircase points. Points with different labels have no common staircase upper bound, so directedness forces one label $x$. If the set of indices is bounded, it has a maximum. If it is unbounded, every member of $E$ is below $\widehat{x}$. Any upper bound must then have the form $\widehat y$, and $y\in U_n(x)$ for every index occurring in $E$. Unboundedness and the decreasing property of the sets $U_n(x)$ give
\[
  y\in\bigcap_{n\in\N}U_n(x)={\downarrow} x.
\]
Thus $\widehat{x}\le\widehat{y}$, proving that $\widehat{x}$ is the supremum.
\end{proof}

Let
\[
  \widehat L=\{\widehat x:x\in L\}\subseteq D_L,
  ~ e(x)=\widehat x.
\]
The point set $\widehat L$ is upper in $D_L$, hence saturated in $\Sigma~\!\!D_L$.

\begin{prop}\label{prop:lower-embedding}
The map
\[
  e:(L,\Lower(L))\longrightarrow \widehat L\subseteq\Sigma~\!\!D_L
\]
is a homeomorphism onto the saturated subspace $\widehat L$.
\end{prop}

\begin{proof}
The map is bijective. Let $O$ be Scott open in $D_L$ and put
\[
  V=e^{-1}(O\cap\widehat L).
\]
If $x\in V$, the directed chain $x_0\le x_1\le\cdots$ has supremum $\widehat x$. Scott inaccessibility gives $x_n\in O$ for some $n$. If $y\in U_n(x)$, then $x_n\le\widehat y$, so $\widehat y\in O$. Hence $U_n(x)\subseteq V$, and $V$ is lower open.

Conversely, let $V$ be lower open and define
\[
  O_V=\{\widehat y:y\in V\}
       \cup\{x_n:U_n(x)\subseteq V\}.
\]
The set $O_V$ is upper. If $\widehat y\in O_V$ and $\widehat y\le\widehat z$, then $z\le_L y$, and lowerness of $V$ gives $z\in V$. If $x_n\in O_V$ and $x_n\le x_m$, then $U_m(x)\subseteq U_n(x)\subseteq V$. If $x_n\in O_V$ and $x_n\le\widehat y$, then $y\in U_n(x)\subseteq V$. These are all possible cases.

Let $E\subseteq D_L$ be directed and suppose $\bigvee E\in O_V$. If $E$ is contained in one staircase and its indices are bounded, its supremum is its maximum and belongs to $E\cap O_V$. If its indices are unbounded, say on the $x$ staircase, then $\bigvee E=\widehat x$ and $x\in V$. Choose $N$ with $U_N(x)\subseteq V$. Some occurring index is at least $N$, so $E\cap O_V\ne\varnothing$.

It remains to consider the case in which $E$ contains point elements. By Lemma~\ref{lem:DL-dcpo},
\[
  \bigvee E=\widehat m,
  ~ m=\bigwedge E^\sharp.
\]
The set $E^\sharp$ is filtered in $L$, hence directed in $L^{\mathrm{op}}$, and its supremum in $L^{\mathrm{op}}$ is $m$. Since $m\in V$ and $V\in\Lower(L)\subseteq\Scott(L^{\mathrm{op}})$ by \eqref{eq:lower-scott-op}, Scott inaccessibility gives some $y\in E^\sharp\cap V$. Thus $\widehat y\in E\cap O_V$.

Hence $O_V$ is Scott open. By construction, $O_V\cap\widehat L=e[V]$, which proves equality of the two topologies.
\end{proof}

\begin{cor}\label{cor:DL-noncosober}
The Scott space $\Sigma~\!\!D_L$ is not co-sober.
\end{cor}

\begin{proof}
By Corollary~\ref{cor:lower-witness} and Proposition~\ref{prop:lower-embedding}, $\widehat C=e[C]$ is compact saturated and $k$-irreducible in the point subspace. Lemma~\ref{lem:basic-transfer}(2) gives the same properties in $\Sigma~\!\!D_L$.

The set $\widehat C$ is not principal. If $\widehat C={\uparrow} d$ and $d=x_n$ is a staircase point, then $d\in{\uparrow} d$ but $d\notin\widehat L$. If $d=\widehat x$, then
\[
  {\uparrow}_{D_L}\widehat{x}=\{\widehat{y}:y\le_L x\},
\]
so $C={\downarrow} x$, contradicting \eqref{eq:C-nonprincipal}.
\end{proof}

\section{The countable counterexample}

\begin{proof}[Proof of Theorem~\ref{thm:countable-intro}]
Let
\[
  Y_c=\Alex(D_L).
\]
The space $Y_c$ is countable, $T_0$, and co-sober by Lemma~\ref{lem:alex-facts}. Moreover,
\[
  \Q(Y_c)=\Fin(D_L),
\]
which is countable because $D_L$ is countable. Corollary~\ref{cor:DL-noncosober} and Corollary~\ref{cor:finite-reflection} imply that
\[
  P_{\Scott}(Y_c)=\Sigma\Fin(D_L)
\]
is not co-sober. Lemma~\ref{lem:alex-facts}(3) shows that the upper Vietoris space $P_S(Y_c)$ is the Alexandrov space of the same Smyth ordered set $\Fin(D_L)$ and is therefore co-sober.
\end{proof}

Hence the upper Vietoris and Scott topologies have different co-sobriety behavior on the same countable Smyth power poset.

\section{Concluding remarks}

Theorem~\ref{thm:main} shows that the Scott power construction does not preserve co-sobriety. The finite transfer uses the point-saturation and union maps from the full compact saturated setting and Heckmann and Keimel's finite Rudin result at the compactness step. The separate staircase construction yields a countable example and, for the same countable Smyth power poset, distinguishes the upper Vietoris topology from the Scott topology.

\end{document}